\documentclass{ifacconf} 

\usepackage[utf8]{inputenc}

\usepackage{natbib} 
\usepackage{lipsum}
\usepackage[T1]{fontenc}
\usepackage{enumitem}
\usepackage{amsmath,amssymb}
\usepackage{mathtools}
\usepackage[ruled,vlined,algo2e]{algorithm2e} %must be loaded before cleveref

\usepackage{xcolor}

\usepackage{bm}

\newcommand{\vect}[1]{\boldsymbol{#1}} % Vector
\newcommand{\m}[1]{\mathbf{#1}} % Matrix

\newcommand{\R}{\mathbb{R}}

\newcommand{\mpi}{+} %Moore-Penrose inverse
 
\newcommand{\argmax}{\operatorname{argmax}} 
 
\newcommand{\pdc}[1]{\mathbb{S}_{+}^{#1}} 
\newcommand{\support}[1]{\mathcal{S}(#1)} 
\newcommand{\Int}[1]{\operatorname{int}(#1)}

\newcommand{\subseteqeq}{\subseteq_0} 
\newcommand{\ellipsoid}[2]{\operatorname{E}(#1,#2)} 
\newcommand{\ball}[2]{\operatorname{B}(#1,#2)}

\newtheorem{theorem}{Theorem}
\newtheorem{lemma}{Lemma}
\newtheorem{corollary}{Corollary}
\newtheorem{proposition}{Proposition}

\newtheorem{remark}{Remark}
\renewcommand{\qed}{$\blacksquare$}

\begin{document}
	\begin{frontmatter}
		\title{An Efficient Method to Verify the Inclusion of Ellipsoids}
		
		\thanks[footnoteinfo]{
		RJ is a FNRS honorary Research Associate. This project has received funding from the European Research Council (ERC) under the European Union's Horizon 2020 research and innovation programme under grant agreement No 864017 - L2C. RJ is also supported by the Innoviris Foundation and the FNRS (Chist-Era Druid-net).
		JC is a FRIA Research Fellow. }

		\author[First]{Julien Calbert} 
		\author[First]{Lucas N. Egidio}
		\author[First]{Raphaël M. Jungers} 
		
		\address[First]{ICTEAM Institute, UCLouvain, Louvain-la-Neuve, Belgium\\ (e-mail: \{julien.calbert, lucas.egidio, raphael.jungers\}@uclouvain.be).}
	\begin{abstract}
	    We present a novel method for deciding whether a given $n$-dimensional ellipsoid contains another one (possibly with a different center). This method consists in constructing a particular concave function and deciding whether it has any value greater than $-1$ in a compact interval that is a subset of $[0,1]$. This can be done efficiently by a bisection algorithm method that is guaranteed to stop in a finite number of iterations. The initialization of the method requires $\mathcal{O}(n^3)$ floating-point operations and evaluating this function and its derivatives requires  $\mathcal{O}(n)$. This can be also generalized to compute the smallest level set of a convex quadratic function containing a finite number of $n$-ellipsoids. In our benchmark with randomly generated ellipsoids, when compared with a classic method based on semidefinite programming, our algorithm performs 27 times faster for ellipsoids of dimension $n=3$ and 2294 times faster for dimension $n=100$. We illustrate the usefulness of this method with a problem in the control theory field.
	\end{abstract}
	
	\end{frontmatter}

\section{Introduction}
For many problems in control and estimation theory, ellipsoidal sets represent a sensible compromise between expressiveness power and numerical tractability. Their simple characterization by a convex quadratic function allows the expression of involved control objectives and constraints as optimization problems that can be handled relatively efficiently by {convex optimization}. Also, they are suitable for characterizing uncertainties and disturbances, particularly when Gaussian noise is assumed. For a non-exhaustive list of classic applications on control, we refer to~\citep{boyd1994linear,kurzhanski1997ellipsoidal} and, for estimation, a few instances are~\citep{schweppe1968recursive,bertsekas1971recursive,zolghadri1996algorithm}.

More recently, with the development of modern control techniques, such as abstraction-based control design, neural-network-based control, and data-driven control, among others, the representation of mathematical concepts through  ellipsoidal sets has been shown to be also useful in these contexts, e.g., developing barrier functions and local controllers~\citep{he2020bp,egidio2022state}, assessing the safety of neural-networks \citep{fazlyab2019probabilistic}, and capturing data uncertainties in data-driven control methods~\citep{bisoffi2022data}.

In view of all these applications and the growing necessity for efficient methods to perform numerical operations with ellipsoids, we present a novel approach to verify whether one $n$-ellipsoid $\mathcal{E}\subset\R^n$ is a subset of another ellipsoid $\mathcal{E}_0\subset\R^n$. Our contributions in this paper are listed as follows:
\begin{itemize}[leftmargin=*]
    \item we write the problem of inclusion of two ellipsoids as a concave minimization problem for which strong duality holds. Then, one can decide the inclusion by computing the maximum of the \emph{dual function}, which is a scalar concave function that can be evaluated in $O(n)$ floating-point operations (FLOPs), as well as its derivative.
    \item we prove that the dual search domain for the maximum of the dual function can be restricted to a compact set contained in the interval $[0,1]$, which makes it suitable to be handled by bisection algorithms. Additionally, an early stop criterion is presented, which is triggered within a finite number of iterations of the bisection algorithm when strict inclusion holds.
    \item we generalize our algorithm to compute the smallest level set of a positive definite quadratic function containing a finite number of $n$-ellipsoids. This is a problem with applications in control theory, and we present an example, namely, calculating control forward invariant sets.
    \item we show that, in a benchmark consisting of randomly generated ellipsoids, when compared to the classic semidefinite programming-based method, our approach performs, on average, about 27 times faster for ellipsoids of dimension $n=3$ and 2294 times faster in dimension $n=100$. 
\end{itemize}

\textbf{Literature Review:}
The most simple case of testing the ellipsoid inclusion happens when they share the same  center. As it will be discussed, this can be solved by comparing the eigenvalues of the Hessian matrices of the quadratic functions defining each of them. In the general case, when the ellipsoids do not share the same center, the inclusion problem can be reformulated (see \cite[p.~43]{boyd1994linear}) as a linear matrix inequality (LMI) problem using the S-lemma~\cite[Thm.~2.2]{polik2007survey}.

This LMI problem yields a semidefinite program (SDP) and, thus, can be solved by one of many readily available SDP solvers. Nonetheless, the performance and accuracy of these solvers are often not ideal as they do not leverage the specific structure of the problem being solved. Therefore, in this work, we design a method that, by exploiting the structure of the ellipsoidal inclusion problem, outperforms general-purpose SDP solvers such as SDPA~\citep{yamashita2010high} and Mosek~\citep{mosek}.

In a similar fashion,  the authors of~\cite[Prop.~2]{gilitschenski2012robust} reformulate the ellipsoid \emph{intersection} problem as the minimization of a convex scalar function in a bounded interval. Their method requires $O(n^4)$ FLOPs and is based on algebraic geometry.
For this same problem of intersection of $n$-ellipsoids, in \citep{ros2002ellipsoidal}, the authors present an algorithm to compute, among all ellipsoids that are convex combinations of two given ones, $\mathcal{E}$ and $\mathcal{E}_0$,  the one with minimal volume and that contains $\mathcal{E}\cap\mathcal{E}_0$. Their method relies on the computation of the root of a polynomial of degree $2n-1$ defined in a bounded interval, see \cite[Thm.~3]{ros2002ellipsoidal}. 

For the problem of \emph{inclusion} of ellipsoids, however, no tailored procedure is available in the literature to the best of our knowledge.

\textbf{Notations:}
We denote a matrix $\m A\in\R^{n\times m}$ by a capital letter and in bold, and the element of the $i$-th row and $j$-th column of $\m A$ (with $1\leq i\leq n$ and $1\leq j \leq m$) by $A_{ij}$.
The Moore–Penrose inverse of a matrix $\m A$ is denoted as $\m A^\mpi$.
 Given a vector $\vect v \in \R^n$, we define the support of vector $\vect v$ as $\support{\vect v} = \{i\in\{1,\ldots,n\} : v_i\ne 0\}$. 
We denote by $\lambda_{\min}(\m A)$ the smallest eigenvalue in absolute value of the matrix $\m A\in\R^{n\times n}$.
We denote by $\pdc{n}$ the set of positive definite matrices of dimension $n$. Also, $\m A\succ \m 0$  represents that  $\m A\in\mathbb{S}_+^n$ and $\m A \succeq \m 0$, that $\m A\in \partial \mathbb{S}_+^n$, i.e., $\m A$ is positive semidefinite.
The convex hull of a set of vectors $\vect v_1,\ldots,\vect v_m\in \R^n$ is denoted ${\rm co}\{\vect v_1,\ldots,\vect v_m\}$. Given a set $S\subset\R^n$, we denote by $\Int{S}$ and $\partial S$, the interior and the boundary of $S$, respectively.
Expressions containing the symbol ``$\pm$'' should be read twice replacing it by ``$+$'' and ``$-$''.

\textbf{Outline:}
This paper is structured as follows. Section \ref{sec:preliminary} includes the mathematical background needed throughout this paper. Section \ref{sec:main-results} is devoted to our main result: the optimization formulation of the inclusion test. In Section \ref{sec:numerical} we provide the details of the practical implementation, benchmarking with off-the-shelf solvers, and an example of application in control theory.

\section{Preliminary Results}\label{sec:preliminary}
% \begin{lemma}
%   Consider two quadric $n$-solids $Q_1\subseteq Q_2$. $S=\{1/\lambda_{\min} \} \iff \partial Q_1 \cup \partial Q_2 \neq \emptyset$  
% \end{lemma}
\subsection{The Inclusion of Ellipsoids}
Before presenting our main results, we first state some definitions and present existing results regarding the problem of verifying the inclusion of ellipsoids. An $n$-ellipsoid with \emph{center} $\vect{c}\in\R^n$ and \emph{shape} defined by $\m P\succ \m 0$ is denoted as 
\begin{equation}
\label{eq:def_ellipsoid}
    \ellipsoid{\vect c}{\m P} \coloneqq  \{\vect x\in\R^n: (\vect x-\vect c)^\top \m P(\vect x-\vect c) \le 1\}.
\end{equation}
Naturally, the $n$-dimensional Euclidean ball of radius $r>0$ and centered at ${\vect c}$ is denoted as $\ball{{\vect c}}{r}\coloneqq\ellipsoid{{\vect c}}{r^{-{1}/{2}}{\m I}_n}$. For two ellipsoids $\mathcal{E} = \ellipsoid{\vect c}{\m P} $ and  $\mathcal{E}_0 = \ellipsoid{\vect c_0}{\m P_0}$, the \emph{inclusion}, the \emph{strict inclusion} and the \emph{non-inclusion} are denoted by the standard mathematical symbols $\subseteq,~\subset,$ and $\not\subseteq$. Besides these, an additional relation $\subseteqeq$ is considered in this paper and defined as follows
\begin{equation}
    \mathcal{E} \subseteqeq \mathcal{E}_0 \iff \mathcal{E} \subseteq \mathcal{E}_0~\text{and}~ \partial \mathcal{E} \cap   \partial \mathcal{E}_0 \neq \emptyset,
\end{equation}
which  means that $\mathcal{E}$ is included in $ \mathcal{E}_0$ and both ellipsoids have common points in their boundaries that we will call \emph{contact points}. For studying the inclusion in our context, this situation denotes an extreme case and yields a particular interpretation of our algorithm to be presented. 

Verifying whether one ellipsoid is included in another can be equivalently rewritten as verifying if a surrogate ellipsoid is included in a unit Euclidean ball centered at the origin. The next lemma formalizes this equivalence.

\begin{lemma}\label{lem:changeVariables}
Let matrices $\m P,\m P_0\in\pdc{n}$ and vectors ${\vect c},~{\vect c}_0 \in\R^n$ be given. The following equivalences hold
\begin{align*}
    \ellipsoid{\vect c}{\m P} \subset \ellipsoid{\vect c_0}{\m P_0} &\Leftrightarrow \ellipsoid{\vect{\tilde{c}}}{\m{\tilde{P}}}\subset \ball{\vect 0}{1};\\
    \ellipsoid{\vect c}{\m P}\subseteqeq \ellipsoid{\vect c_0}{\m P_0}&\Leftrightarrow \ellipsoid{\vect{\tilde{c}}}{\m{\tilde{P}}}\subseteqeq \ball{\vect 0}{1};
\end{align*}
with 
\begin{equation}\label{eq:tildes}
    \vect{\tilde{c}} = \m L_0^\top (\vect c-\vect c_0),\quad \m{\tilde{P}} = \m L_0^{-1}\m P \m L_0^{-\top}
\end{equation}
and $\m L_0$ defines the Cholesky factorization of $\m P_0 = \m L_0 \m L_0^\top$.
\end{lemma}
\begin{pf}
As $\m P_0 \succ \m 0$, we have that $\m L_0$ is regular. By applying the change of variables $\vect{\tilde{x}} = \m L_0^\top (\vect x-\vect c_0)$, we have that $\ellipsoid{\vect{c}}{\m{P}}$ and $\ellipsoid{\vect{c}_0}{\m{P}_0}$ becomes respectively $\ellipsoid{\vect{\tilde{c}}}{\m{\tilde{P}}}$ and $\ball{\vect 0}{1}$ in the space of $\tilde{\vect x}$.\hfill\hfill\qed
\end{pf}
As a consequence, in this paper, we will equivalently study the problem of verifying if an ellipsoid is included in a Euclidean $n$-ball of radius $1$ given that an appropriate change of variables transforming the original problem into this one always exists. Notice that this can be done under $\mathcal{O}(n^3)$ arithmetic operations because of the required Cholesky factorization~\citep{higham2009cholesky}.

For the sake of completeness, before presenting the necessary and sufficient condition for inclusion on which our method is based, we will discuss other simpler criteria that allow us to sufficiently determine whether one ellipsoid is included in the other or not. These can be used as preliminary tests before running our algorithm to further speed up the execution time of an inclusion verification routine.
\begin{proposition}\label{prop:pretests}
Let matrices $\m P, \m P_0\in\pdc{n}$ and vectors ${\vect c},~{\vect c}_0 \in\R^n$ be given. 
The following are necessary conditions for $\ellipsoid{\vect c}{\m P} \subseteq \ellipsoid{\vect c_0}{\m P_0}$:
\begin{enumerate}
    \item $\vect{c}\in \ellipsoid{\vect c_0}{\m P_0}$,
    \item $\m P\succeq \m P_0$.
\end{enumerate}
Moreover, if $\vect c = \vect c_0$ then condition (2) is also sufficient.
\end{proposition}
\begin{pf}
Let us demonstrate each of these two statements.
\begin{enumerate}
    \item This trivially holds from the fact that $\vect c\in \ellipsoid{\vect c}{\m P}$.
    \item To show a contradiction, assume that $\ellipsoid{\vect c}{\m P} \subseteq \ellipsoid{\vect c_0}{\m P_0}$ but also that there exists $\vect{v}\in\R^n$ such that $\vect{v}^\top\m P_0\vect{v}> \vect{v}^\top\m P\vect{v}=1$, without loss of generality. Therefore, $\vect x_+, \vect x_- \in\ellipsoid{\vect c}{\m P}$, with $\vect x_+=\vect c+\vect v,~\vect x_-=\vect c-\vect v$. On the other hand, $\vect x_+$ and/or $\vect x_-$ are not in $\ellipsoid{\vect c_0}{\m P_0}$, which can be verified by developing the left-hand side expression in the definition~\eqref{eq:def_ellipsoid} as
    \begin{align}
        (\vect{x}_\pm \!-\!\vect{c}_0)^\top\m P_0(\vect{x}_\pm \!-\!\vect{c}_0) =&(\vect{c} \!-\!\vect{c}_0\!\pm\!\vect v)^\top\m P_0(\vect{c} \!-\!\vect{c}_0\!\pm\!\vect v) \nonumber\\
        >&(\vect{c} \!-\!\vect{c}_0)^\top \m P_0 (\vect{c} \!-\!\vect{c}_0) \nonumber\\
        &~~\pm 2\vect{v}^\top\m P_0 (\vect{c} -\vect{c}_0)  + 1  \nonumber \\
        >&\pm 2\vect{v}^\top\m P_0 (\vect{c} -\vect{c}_0)  + 1.
    \end{align}
    This shows that at least $\vect x_+\notin \ellipsoid{\vect c_0}{\m P_0}$ or $\vect x_-\notin \ellipsoid{\vect c_0}{\m P_0}$, which is a contradiction. 
\end{enumerate}
To show the sufficiency of (2) for $\vect c=\vect c_0$, note that
    $(\vect x - \vect c)^\top \m P (\vect x - \vect c)\geq (\vect x - \vect c)^\top \m P_0 (\vect x - \vect c) $ for all $\vect x\in\R^n$. Therefore, for all $\vect x\in\ellipsoid{\vect c}{\m P}$ we have $1\geq(\vect x - \vect c)^\top \m P (\vect x - \vect c)\geq (\vect x - \vect c)^\top \m P_0 (\vect x - \vect c)$, which implies that $\vect x \in\ellipsoid{\vect c_0}{\m P_0}$.
\hfill\hfill\qed
\end{pf}
Although simple, these tests allow us to efficiently decide about the inclusion of ellipsoids for some cases. Additionally, for the case $\vect c = \vect c_0$, one can show in a similar fashion that $\m P \succ \m P_0$ implies strict inclusion.

\subsection{An optimization approach }

Let us introduce the following optimization problem 
 \begin{equation}\label{eq:primal}
\begin{aligned}
p^*  = &\min_{\vect x\in \R^n} \ -\vect x^\top\vect x\\
&\textrm{s.t.} \quad  (\vect x-\vect c)^\top \m P (\vect x-\vect c) \le 1,
\end{aligned}
\end{equation}
which finds the point $\vect{x}$ of maximum Euclidean norm within $\ellipsoid{\vect c}{\m P}$. Therefore, the optimal value $p^*$ is related to the problem of verifying the inclusion of an ellipsoid inside the euclidean unit ball as follows:
\begin{align}
    \ellipsoid{\vect c}{\m P} \subseteq \ball{\vect 0}{ 1} &\Leftrightarrow p^* \ge -1.
\end{align} 
It is also straightforward to show that $p^*=-1$ if and only if $\ellipsoid{\vect c}{\m P} \subseteqeq \ball{\vect 0}{1}$.
Although~\eqref{eq:primal} is a non-convex optimization problem, it has some noteworthy properties. First, there always exists a (strictly) feasible point, given that the ellipsoid $\ellipsoid{\vect{c}}{\m P}$ contains its center in its interior. Also, it is always bounded, given that it is a minimization of a concave function within a convex set~\citep[p.~342]{rockafellar1970convex}. Finally, the \emph{Lagrangian} of this optimization problem is a quadratic function on $\vect{x}$, given as 
\begin{equation}\label{eq:lagrangian}
\mathcal{L}(\vect x,\beta) = -\vect x^\top \vect x  + \beta \left( (\vect x-\vect c)^\top \m P (\vect x-\vect c)-1\right)
\end{equation}
where $\beta$ is the Lagrange multiplier associated with the unique constraint of \eqref{eq:primal}.
Naturally, the \emph{Lagrange dual function} $g:\R^n\times \R\rightarrow \R$ is defined by
\begin{equation}\label{eq:dual_function}
    g(\beta) = \min_{\vect x\in \R^n} \ \mathcal{L}(\vect x,\beta)
\end{equation}
and the dual optimization problem by
\begin{equation}\label{eq:dual}
    d^* = \max_{\beta\in \mathcal{D}_g} g(\beta),
\end{equation}
where the dual domain is defined as
\begin{equation}
	\mathcal{D}_g =\{\beta\geq 0~:~ g(\beta) >-\infty\}.
\end{equation}

Notice that, because of its quadratic nature, the lower-boundedness of the Lagrangian function is closely related to the sign of its Hessian 
\begin{equation}
	\nabla_{\vect x}^2\mathcal{L}(\vect x,\beta) = \beta\m P-\m I .
\end{equation}
As discussed in~\cite[p.~458]{boyd2004convex}, for a given $\beta\geq 0$, this function is bounded from below if $\nabla_{\vect x}^2\mathcal{L}(\vect x,\beta)  \succ \m 0$
or if $\nabla_{\vect x}^2\mathcal{L}(\vect x,\beta) \succeq \m 0 $ and there exists $\vect{x}\in\R^n$ such that $\nabla_{\vect{x}}\mathcal{L}(\vect x,\beta)=\vect 0$. This implies that the domain of the dual function is either $\mathcal{D}_g=[1/\lambda_{\min}(\m P),\infty)$ or $\mathcal{D}_g=(1/\lambda_{\min}(\m P),\infty)$, depending on the matrix $\m P$ and the vector $\vect{c}$ defining the Lagrangian~\eqref{eq:lagrangian}. 

The next section clarifies how the dual function~\eqref{eq:dual_function} can be used to construct an efficient algorithm for verifying the inclusion of ellipsoids.

\section{Main results}\label{sec:main-results}
\subsection{An Algorithm to Test the Inclusion of Ellipsoids}

 Corresponding to the problem of verifying  whether the inclusion $\ellipsoid{\vect{c}}{\m{P}}\subseteq\ball{\vect{0}}{1}$ holds, we define the infinitely differentiable function $\ell_{\vect c,\m P}:\mathcal{I}_{\m P}\rightarrow \R$ as 
\begin{equation}\label{eq:function_l}
    \ell_{\vect c,\m P}(\beta) \coloneqq  -\beta - \sum_{i\in \support{\vect{\bar{c}}}} \bar{c}_i^2 \frac{\lambda_i \beta}{\lambda_i \beta-1} 
\end{equation}
and its domain 
\begin{equation}\label{eq:Ip}
	\mathcal{I}_{\m P}\!\coloneqq\!\!\left\{ \begin{array}{rl}
	\!\!(\lambda_{\min}(\m P)^{-1},\infty),& ~\text{if} ~\exists i\! \in\!\mathcal{S}(\bar{\vect{c}}),~\lambda_i \!=\! \lambda_{\min}(\m P)\\ 
\!\!\!\ [\lambda_{\min}(\m P)^{-1},\infty),& ~\text{otherwise}
	\end{array}\right.
\end{equation}
where $\vect{\bar{c}} = \m V^\top \vect c$, $\lambda_i = D_{ii}$ and $(\m V, \m D)$ constitute the spectral decomposition of $\m P = \m V \m D \m V^\top$. Notice that, within the domain of $\ell_{\vect c,\m P}(\beta)$, the only case where the expression~\eqref{eq:function_l} is undefined is when $\beta=1/\lambda_{\min}(\m P)$ and there exist $i\in\support{\vect{\bar{c}}}$ such that $\lambda_i=\lambda_{\min}(\m P)$. In this case, we remove this lower limit from its domain. In the following lemma, we present some important properties of $\ell_{\vect c,\m P}(\beta)$.

\begin{lemma}\label{lemma:properties}
Let $\m P\in \pdc{n}$ and $\vect{c}\in\R^n$ be given. The function $\ell_{\vect c,\m P}(\beta)$, given in~\eqref{eq:function_l}, has the following properties:
\begin{enumerate}
    \item $\ell_{\vect c,\m P}(\beta)$ is concave in its domain $\mathcal{I}_{\m P}$;
    \item $\forall \beta\in\mathcal{I}_{\m P},~ \ell_{\vect c,\m P}(\beta)<0$;

\end{enumerate} 
\end{lemma}
\begin{pf}

For $\beta\in \mathcal{I}_{\m P}$, we decompose $\ell_{\vect c,\m P}$ as 
$$\ell_{\vect c,\m P}(\beta)= h_0(\beta) + \sum_{i\in\support{\vect{\bar{c}}}} \ \bar{c}_i^2 \lambda_{i} h_i(\beta)$$
with $h_0(\beta) = -\beta$ and $h_i(\beta) = {-\beta}/({\lambda_i \beta}-1)$,~$i\in\support{\vect{\bar{c}}}$. Note that, $\bar{c}_i^2 \lambda_i\ge0$  for all $i$, given that $\m P\succ \m 0$. Below we proof each property in the statement.
\begin{enumerate}
\item Evaluating the first and second derivatives of $\ell_{\vect c,\m P}(\beta)$ for $\beta\in\mathcal{D}_g$, we have
\begin{align}
    \ell_{\vect c,\m P}'(\beta) &= -1+ \sum_{i\in\support{\vect{\bar{c}}}} \bar{c}_i^2 \frac{\lambda_i}{(\lambda_i\beta-1)^2},\label{eq:dl}\\
    \ell_{\vect c,\m P}''(\beta) &= -2\sum_{i\in\support{\vect{\bar{c}}}} \bar{c}_i^2 \frac{\lambda_i^2}{(\lambda_i\beta-1)^3}.\label{eq:ddl}
\end{align}
Since $\m P\succ \m 0$, we have  $\lambda_i>0$, which implies that the function $\ell_{\vect c,\m P}''(\beta)$ is strictly negative for all $\beta>\lambda_{\min}(\m P)^{-1}$. Thus, $\ell_{\vect c,\m P}(\beta)$ is concave in $\mathcal{I}_{\m P}$.

\item The function $h_0(\beta)<0$ in $\mathcal{I}_{\m P}$ and,  for all $i\in \support{\vect{\bar{c}}}$, the function $h_i(\beta)< 0$ in $({\lambda_i^{-1}},\infty)\supseteq \mathcal{I}_{\m P}$.
We conclude that $\ell_{\vect c,\m P}(\beta)< 0$ in $\mathcal{I}_{\m P}$.
\end{enumerate}
The proof is concluded.\hfill\hfill\qed
\end{pf}

As demonstrated in the previous lemma, inside its domain $\mathcal{I}_{\m P}$ this function is concave. Hence,  we can obtain 
\begin{equation}
    \ell_{\vect c,\m P}^* \coloneqq \sup_{\beta \in \mathcal{I}_{\m P}}\ \ell_{\vect c,\m P}(\beta)\label{eq:max_function_l}
\end{equation}
by classic optimization algorithms such as Newton's method, which guarantees a quadratic convergence rate to $\ell_{\vect c,\m P}^*$ \cite[Thm.~3.5]{nocedal1999numerical}, or by a bisection algorithm, which avoids the computation of the second derivative. The following theorem connects the function~\eqref{eq:function_l} with the ellipsoidal inclusion problem.

\begin{theorem}\label{th:ellipsidInclusionUnitBall}
Let an ellipsoid $\ellipsoid{\vect{c}}{\m{P}}$ be given. Define the function $\ell_{\vect c,\m P}(\beta)$ as in~\eqref{eq:function_l} and consider its supremum $\ell_{\vect c,\m P}^*$ over the domain $\mathcal{I}_{\m P}$. The following equivalences always hold:
\begin{align*}
    \ellipsoid{\vect c}{\m P}\subset\Int{\ball{\vect 0}{1}} \Leftrightarrow&\ \ell_{\vect{c},\m{P}}^* > -1;\\
    \ellipsoid{\vect c}{\m P}\subseteqeq\ball{\vect 0}{1} \Leftrightarrow&\ \ell_{\vect{c},\m{P}}^* = -1;
\end{align*}
\end{theorem}

\begin{pf} First, let us recall that the optimal solution $p^*$ of the primal optimization problem~\eqref{eq:primal} allows us to decide whether $\ellipsoid{\vect c}{\m P}$  is included or not inside $\ball{\vect 0}{1}$, i.e., the inclusion $ \ellipsoid{\vect c}{\m P}\subset\Int{\ball{\vect 0}{1}}$ (resp. $\ellipsoid{\vect c}{\m P}\subseteqeq\ball{\vect 0}{1}$) holds if and only if $p^*>-1$ (resp. $p^*=-1$). Let us show that this problem has no duality gap, that is, $p^*=d^*$ where $d^*$ is the optimal solution of the dual problem~\eqref{eq:dual}. Notice that Slater’s constraint qualification~\cite[p.~226]{boyd2004convex} holds  given that $\vect x = \vect c$ is a strictly feasible point. Although the primal problem is not convex due to the concave objective function, Slater's condition implies strong duality in this case, given that this problem is the minimization of a quadratic function subject to a single quadratic inequality, see \cite[Sec.~B4]{boyd2004convex} for a detailed proof. 

Therefore, we must now show that $\ell_{\vect c,\m P}^*=d^*$. Notice that, the minimization problem within the definition of the dual function~\eqref{eq:dual_function} is convex for $\beta\in\mathcal{D}_g$ and, therefore, its global minimizers $\vect{x}^*$ for a given $\beta\in\mathcal{D}_g$ are all points satisfying the first-order optimality condition 
\begin{equation}\label{eq:suff_cond_min_lagrangian}
    \nabla_{\vect x} \mathcal{L}(\vect x^*,\beta) = 2(\beta \m P-\m I) \vect{x}^*  + 2\beta \m P \vect c = \vect 0.
\end{equation}
If $\beta>1/\lambda_{\min}(\m{P})$, then the matrix $(\beta \m P-\m I)$ can be inverted and there is a unique minimizer $\vect x^*(\beta) = -(\beta \m P- \m I)^{-1} \beta \m P \vect c$. However, when $\beta=1/\lambda_{\min}(\m{P})$, two cases may occur: 1) there is no $\vect{x}^*$ such that \eqref{eq:suff_cond_min_lagrangian} holds and therefore, $1/\lambda_{\min}(\m{P})$ does not belong to the domain $\mathcal{D}_g$ of the dual function $g(\beta)$; 2) there exist infinitely many $\vect{x}^*$ satisfying~\eqref{eq:suff_cond_min_lagrangian}, which are given by
\begin{equation}\label{eq:mpiSolution}
    \vect x^*(\beta) = -(\beta \m P- \m I)^\mpi \beta \m P \vect c + \m{N}\vect{v}
\end{equation}
where $\m{N}\in\R^{n\times m}$ has columns spanning the $m$-dimensional nullspace of $\big(\m P-\lambda_{\min}(\m{P})\m I\big)$ and $\vect{v}\in\R^m$. As each of these minimizers is a global minimizer, without loss of generality, let us choose $\vect x^*(\beta) = -(\beta \m P- \m I)^\mpi \beta \m P \vect c$ (i.e., $\vect{v}=\vect{0}$). Naturally, this minimizer must satisfy~\eqref{eq:suff_cond_min_lagrangian} by assumption. Substituting $\vect x^*(\beta)$ into~\eqref{eq:suff_cond_min_lagrangian} and performing a few algebraic manipulations yields 
\begin{align}
	-\beta(\beta \m P-\m I)^\mpi  \m P \vect c  &= \beta \m P \vect c-\beta^2 \m P (\beta \m P-\m I)^\mpi \m P \vect c .\label{eq:condition_eq_0}
\end{align}
 Hence, the dual function~\eqref{eq:dual_function} can be rewritten as
\begin{align*}		
	g(\beta) &=  \mathcal{L}(\vect x^*(\beta),\beta)\\
	&=\beta (\vect c^\top \m P\vect c-1) - \beta^2 \vect c^\top \m P (\beta \m P-\m I)^{\mpi} \m P \vect c\\
	&=-\beta  + \vect c^\top (\beta \m P - \beta^2 \m P(\beta \m P-\m I)^{\mpi} \m P) \vect c \\
	&=-\beta  - \beta\vect c^\top  (\beta \m P-\m I)^{\mpi} \m P \vect c \\	
	&=-\beta  - \beta\vect{\bar{c}}^\top \m  (\beta \m D-\m I)^{\mpi}\m{D} \vect{\bar{c}},
\end{align*}
where~\eqref{eq:condition_eq_0} was used to obtain the before last equation and the last equation uses the spectral decomposition $\m P=\m V \m D \m V^\top$ and the transformation $\vect{\bar{c}} = \m V^\top \vect c$. Moreover, due to the fact that all matrices in this last expression are diagonal, one can rewrite
\begin{equation}
	g(\beta) = -\beta -\sum_{i\in \support{\vect{\bar{c}}}}\bar{c}_i^2 \frac{\lambda_i\beta}{\lambda_i\beta-1}.
\end{equation}
Therefore the dual function $g(\beta)=\ell_{\vect c,\m P}(\beta)$ for all $\beta\in \mathcal{D}_g=\mathcal{I}_{\m P}$, and, hence, $d^* =\ell_{\vect c,\m P}^*$, concluding the proof.
\hfill\hfill\qed
\end{pf}

\begin{figure}
    \centering
    \includegraphics[width=\linewidth]{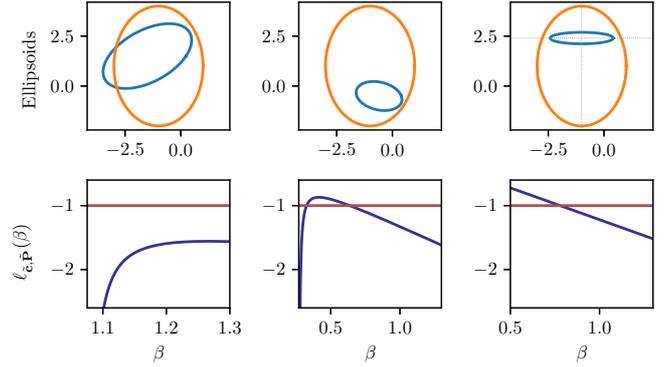}
    \caption{Illustration of results of Theorem~\ref{th:ellipsidInclusionUnitBall}. The ellipsoid $\mathcal{E}$ (blue) is contained in $\mathcal{E}_0$ (orange) if and only if the maximum of $\ell_{\tilde{\vect{c}},\tilde{\m{P}}}(\beta)$ is greater than $-1$.} 
    \label{fig:ellipsoids_and_ell}
\end{figure}

\begin{remark}
The inclusion $\ellipsoid{\vect c}{\m P}\! \subseteq\! \ball{\vect 0}{1}$ is equivalent to 
\begin{equation}
    \min_{\beta\ge 0} \ \max_{\vect x\in \R^n} \mathcal{L}(\vect x,\beta) \ge -1
\end{equation}
which, in turn, is equivalent to 
\begin{equation}\label{eq:geometric_interpretation}
\exists\beta\geq 0,~\forall \vect x\in \R^n :\ \mathcal{L}(\vect x,\beta) \ge -1.
\end{equation}

We can rewrite $\mathcal{L}(\vect x,\beta) + 1 =  [\vect x^\top ~~1] \m F(\beta)[\vect x^\top ~~1]^\top$ with 
\begin{equation}\label{eq:A_matrix}
     {\m F(\beta)}\coloneqq \begin{pmatrix}
    \beta \m P -\m I & -\beta \m P\m c\\
    -\beta \vect c^\top \m P & \beta (\vect c^\top \m P\vect c - 1) +1
     \end{pmatrix}.
\end{equation}
to show that \eqref{eq:geometric_interpretation} is equivalent to
\begin{equation}
    \exists \beta \ge 0:\  \m F(\beta) \succeq \m 0.
\end{equation}
Note that this last inequality is the LMI condition proposed by \citep[Sec.~3.7.1]{boyd1994linear} characterizing the inclusion $\ellipsoid{\vect c}{\m P} \subseteq \ball{\vect 0}{1}$. This shows the equivalence and the connection between this approach and ours. A numerical comparison, in terms of computational time and memory required to verify inclusions by both methods, is presented in Section~\ref{sec:numerical}.
\hfill\hfill $\triangle$
\end{remark}

Theorem~\ref{th:ellipsidInclusionUnitBall} is the foundation of our algorithm to verify the inclusion of ellipsoids. We recall that, although one ellipsoid is considered to be a Euclidean $n$-ball $\ball{\vect{0}}{1}$, Lemma~\ref{lem:changeVariables} provides a change of variables that always allows us to transform the general problem into the one tackled in Theorem~\ref{th:ellipsidInclusionUnitBall}. 

In Figure~\ref{fig:ellipsoids_and_ell}, an illustration of the results of Theorem~\ref{th:ellipsidInclusionUnitBall} is provided. There, three cases of ellipsoids $\mathcal{E}=\ellipsoid{\vect{c}}{\m P}$ (blue) and $\mathcal{E}_0=\ellipsoid{\vect{c}_0}{\m P_0}$ (orange) are depicted, along with the corresponding functions $\ell_{\tilde{\vect{c}},\tilde{\m{P}}}$ in the interval $[1/\lambda_{\min}(\tilde{\m{P}}), 1.3]$, where $\tilde{\vect{c}}$ and $\tilde{\m{P}}$ are given in~\eqref{eq:tildes}. Notice that, for the first case (left) the inclusion of $\mathcal{E}$ (blue ellipsoid) within $\mathcal{E}_0$ (orange) does not hold and the corresponding function $\ell_{\tilde{\vect{c}},\tilde{\m{P}}}$ is always strictly below $-1$. For the second case (middle), the inclusion holds, and therefore, the maximum of $\ell_{\tilde{\vect{c}},\tilde{\m{P}}}$ is greater than $-1$. Finally, a third case (right) illustrates the case when $1/\lambda_{\min}(\tilde{\m{P}}) \in\mathcal{D}_g$. This happens  because, after transforming $\mathcal{E}_0$ into the unit Euclidean ball centered at the origin, the center $\tilde{\vect{c}}$ of the transformed $\mathcal{E}$ is perpendicular to its greatest semi-axis (or semi-axes), which is the eigenvector associated to $\lambda_{\min}(\tilde{\m{P}})$.

Before introducing our general algorithm, let us present an additional property of the scalar function $\ell_{\vect{c}, \m P}(\beta)$, defined in~\eqref{eq:function_l}, which will be important for implementation purposes. 

\begin{proposition}\label{prop:upper_bound}
    Let $\m P \in\pdc{n}$ and $\vect{c}\in\R^n$ be given. The function $\ell_{\vect c,\m P}(\beta)$, defined in~\eqref{eq:function_l} satisfies
 $$ \ell_{\vect c,\m P}(\beta) < -1,\quad \forall \beta>\max\big\{\lambda_{\min}(\m P)^{-1}, 1-\vect c^\top \vect c\big\}.$$
\end{proposition}
\begin{pf}
By contradiction, assume that there exists $\beta_0>1-\vect c^\top \vect c$ such that  $\beta_0> \lambda_{\min}(\m P)^{-1}$ and $\ell_{\vect c,\m P}(\beta_0)\geq -1$. From the proof of Theorem~\ref{th:ellipsidInclusionUnitBall}, we have that $\ell_{\vect c,\m P}(\beta_0) = g(\beta_0)$, which implies \begin{align*}
    \ell_{\vect c,\m P}(\beta_0) &= \min_{\vect x\in \mathcal{I}_{\m P}} \mathcal{L}(\vect x,\beta_0)\nonumber\\ 
    &\le \mathcal{L}(\vect c,\beta_0)\nonumber\\
    &=-\vect{c}^\top \vect{c}-\beta_0 < -1
\end{align*}
which is a contradiction.\hfill\hfill\qed\end{pf}
Besides providing a useful upper bound on the interval on which $\ell_{\vect c,\m P}(\beta)\geq-1$, reducing the search space for $\beta$ that maximizes $\ell_{\vect c,\m P}(\beta)$, Proposition~\ref{prop:upper_bound} also provides a sufficient condition for $\ellipsoid{\vect{c}}{\m {P}} \not\subseteq\ball{\vect 0}{1}$. Indeed, if $1-\vect c^\top \vect c<1/\lambda_{\min}(\m P)$, then $\ell_{\vect c,\m P}(\beta)< -1$ for all $\beta\geq 1/\lambda_{\min}(\m P)$ and, therefore, the inclusion does not hold.

\SetNlSty{textbf}{}{ }\IncMargin{2em} 
\begin{algorithm2e}[b]
  \SetAlgoLined
   \nl \label{line:pretests1}\uIf{\normalfont\texttt{isPreTestConclusive()}}{ 
   \nl \Return  \texttt{PreTestConclusion()}\;}\label{line:pretests2}
   \nl $\m L_0\gets$ {Cholesky factorization of} $\m P_0 = \m L_0\m L_0^\top$\;\label{line:chol}
  \nl $(\vect{\tilde{c}},~\m{\tilde{P}}) \gets \big( \m L_0^\top (\vect c-\vect c_0),\ \m L_0^{-1}\m P \m L_0^{-\top}\big)$\;
   \nl $\m V,\m D\gets$ {Spectral decomposition of} $\m{\tilde{P}} = \m V\m D\m V^\top$\;\label{line:spec}
   \nl $\ell^\ast_{\tilde{\vect c},\tilde{\m P}} \gets \max_{\beta \in \mathcal{I}_{\vect{\tilde{c}},\m{\tilde{P}}}} \ \ell_{\vect{\tilde{c}},\m{\tilde{P}}}(\beta)$\;\label{line:maximization}
   \nl \uIf{$\ell^\ast_{\tilde{\vect c},\tilde{\m P}}>-1$}{
   \nl   \Return{$\mathcal{E} \subset \Int{\mathcal{E}_0}$}\;
   }\nl\uElseIf{$\ell^\ast_{\tilde{\vect c},\tilde{\m P}}=-1$}{
     \nl \Return{$\mathcal{E} \subseteqeq \mathcal{E}_0$}\;
   }
   \nl\uElse{
  \nl\Return{$\mathcal{E} \not\subseteq \mathcal{E}_0$}\;
   }
  \caption{
    Test the inclusion of an ellipsoid $\mathcal{E} = \ellipsoid{\vect c}{\m P}$ in another ellipsoid $\mathcal{E}_0 = \ellipsoid{\vect c_0}{\m P_0}$.%
  }
  \label{algo:inclusion} 
\end{algorithm2e}\DecMargin{2em}

Finally, these results can be joined together into Algorithm~\ref{algo:inclusion}. It is important to highlight that this algorithm starts at lines~\ref{line:pretests1}-\ref{line:pretests2} by verifying the tests of Propositions~\ref{prop:pretests} and~\ref{prop:upper_bound} through the function \texttt{isPreTestConclusive()}. Whenever these tests are not conclusive, the algorithm proceeds to evaluate the necessary and sufficient condition from Theorem~\ref{th:ellipsidInclusionUnitBall}. To do so, notice that the lines~\ref{line:chol}-\ref{line:spec} comprise the initialization of the algorithm and can be computed within $\mathcal{O}(n^3)$ FLOPs, due to the Cholesky Factorization~\citep{higham2009cholesky} and the spectral decomposition~\citep{banks2022pseudospectral}. Line~\ref{line:maximization} consists in the maximization of a concave scalar function $\ell_{\vect{\tilde{c}},\m{\tilde{P}}}(\beta)$, defined in~\eqref{eq:function_l}, on the interval 
\begin{equation}\label{eq:interval_I}
	\mathcal{I}_{\vect{\tilde{c}},\m{\tilde{P}}}\coloneqq [\lambda_{\min}(\tilde{\m P})^{-1}, 1-\vect{\tilde c}^\top\vect{\tilde c}].
\end{equation}

The upper bound of the interval $\mathcal{I}_{\vect{\tilde{c}},\m{\tilde{P}}}$ is determined from Proposition~\ref{prop:upper_bound}. This maximization problem can be solved without difficulty by a bisection algorithm or, more efficiently, by Newton's method. Notice that, for the former, one needs to compute the first and second derivatives of $\ell_{\vect{\tilde{c}},\m{\tilde{P}}}(\beta)$, defined respectively in~\eqref{eq:dl} and~\eqref{eq:ddl}, which can be done in $\mathcal{O}(n)$ FLOPs.

Algorithm~\ref{algo:inclusion} can be early-stopped whenever a $\beta\in \mathcal{I}_{\vect{\tilde{c}},\m{\tilde{P}}}$ is found such that $\ell_{\vect{\tilde{c}},\m{\tilde{P}}}(\beta)>-1$. However, as discussed in the next subsection, computing this maximum is an efficient way to obtain a  distance between the boundaries of the two ellipsoids (when the inclusion holds) or by how much $\mathcal{E}_0$ must be inflated so it contains $\mathcal{E}$. 

\subsection{Consequences of Theorem~\ref{th:ellipsidInclusionUnitBall}}

The first consequence of Theorem~\ref{th:ellipsidInclusionUnitBall} that we discuss in this paper is the fact that whenever $\ellipsoid{\vect c}{\m P}\subseteqeq\ball{\vect 0}{1}$, knowing $\beta>0$ that maximizes $\ell_{\vect c,\m P}(\beta)$ allows us to fully characterize contact points between  $\ellipsoid{\vect c}{\m P}$ and $\partial\ball{\vect 0}{1}$.

\begin{corollary}\label{cor:contact_point}

If $\ellipsoid{\vect c}{\m P}\subseteqeq \ball{\vect 0}{1}$ then, all contact points
\begin{equation}
    \vect{\bar{x}} \in 
\{\vect{x}~:~(\vect{{x}}-\vect c)^\top \m P (\vect{{x}}-\vect c) =\vect{{x}}^\top \vect{{x}} = 1
\}\label{eq:set_of_contact}
\end{equation}
satisfy
\begin{equation} \label{eq:contact-points}
\vect{\bar{x}} = -(\beta^* \m P- \m I)^\mpi \beta^* \m P \vect c + \m{N}\vect{v}
\end{equation}
for some $\vect{v}\in\R^m$
where $\m{N}\in\R^{n\times m}$ has columns spanning the $m$-dimensional nullspace of $\big(\m P-\beta^*\m I\big)$ and
$$\beta^* = \arg\max_{\beta \in \mathcal I_{\vect c,\m P}}\ \ell_{\vect c,\m P}(\beta).$$
\end{corollary}
\begin{pf}
By the fact that strong duality holds, the contact points $\vect{\bar{x}}$ are optimal solutions of the primal problem~\eqref{eq:primal} and, thus, satisfy the equation \eqref{eq:suff_cond_min_lagrangian} for the optimal $\beta^*$ of the dual problem and, thus, satisfy~\eqref{eq:contact-points}.
\hfill\hfill\qed
\end{pf}

Corollary~\ref{cor:contact_point} provides a complete characterization of the contact points between the boundary of the two ellipsoids. Notice that, there exist a unique contact point whenever $\beta^*>1/\lambda_{\min}(\m P)$, which ensures that the matrix $\big(\m P-\beta^*\m I\big)$ has full rank. However, if $\beta^*=1/\lambda_{\min}(\m P)$, the uniqueness is no longer guaranteed. In this case, a contact point can be found by selecting $\vect{v} = \alpha\vect v_0$ with $\vect v_0\in\R^{m}$ and finding the scalar $\alpha$ that solves the quadratic equation generated by evaluating $\vect{\bar{x}}^\top \vect{\bar{x}} = 1$.

As shown in Theorem~\ref{th:ellipsidInclusionUnitBall}, the  maximal value $\ell^*_{\vect c,\m P}$ being greater or lesser than $-1$ provides some geometrical insights regarding the inclusion. Besides that, its magnitude yields additional information, as the next corollary shows.

\begin{corollary}\label{cor:touching}
For any $\vect c,\vect c_0\in \R^n$ and $\m P,\m P_0\in\pdc{n}$, let $\tilde{\vect c}$ and $\tilde{\m P}$ be defined as in~\eqref{eq:tildes} and $\ell_{\tilde{\vect{c}},\tilde{\m{P}}}^*$, defined in~\eqref{eq:max_function_l}. The following inclusions hold
\begin{equation}
    \ellipsoid{\vect c}{\m P}\subseteqeq \ellipsoid{\vect c_0}{(-\ell_{\tilde{\vect{c}},\tilde{\m{P}}}^*)^{-1} \m P_0}\label{eq:consequence_inclusion_1}
\end{equation}
\begin{equation}
    \ellipsoid{\vect d}{ -\ell_{\tilde{\vect{c}},\tilde{\m{P}}}^*\m P}\subseteqeq \ellipsoid{\vect c_0}{\m P_0}\label{eq:consequence_inclusion_2}
\end{equation}
with  $\vect{d} = (-\ell_{\tilde{\vect{c}},\tilde{\m{P}}}^*)^{-1/2} (\vect c-\vect c_0)+\vect c_0$.
\end{corollary}

\begin{pf}
By Lemma~\ref{lem:changeVariables}, we have that
$\ellipsoid{\vect c}{\m P}\subseteqeq \ellipsoid{\vect c_0}{{\gamma^{-1}} \m P_0}$ if and only if $ \ellipsoid{\gamma^{-1/2}\vect{\tilde{c}}}{\gamma\m{\tilde{P}}} \subseteqeq \ball{\vect 0}{1}$
with $\vect{\tilde{c}},\ \m{\tilde{P}}$ defined in~\eqref{eq:tildes}. By its definition in~\eqref{eq:function_l} we have
\begin{equation}
    \ell_{\tilde{\vect{c}},\tilde{\m{P}}}^*= -\beta^*-\sum_{i\in\support{\vect{\bar{c}}}} \bar{c}^2_i\frac{\lambda_i\beta^*}{\lambda_i\beta^*-1}.\label{eq:proof_lstar}
\end{equation}
Let $\beta_\gamma^* = \gamma^{-1}\beta^*$ and $\gamma=-\ell_{\tilde{\vect{c}},\tilde{\m{P}}}^*$. By the property (2) in Lemma~\ref{lemma:properties}, we have $\beta_\gamma^*>0$. Tedious but simple algebraic manipulations where~\eqref{eq:proof_lstar} is used, yield $\ell_{\gamma^{-1/2}\vect{\tilde{c}},\gamma\m {\tilde{P}}}(\beta_\gamma^*)=-1$ and $\ell_{\gamma^{-1/2}\vect{\tilde{c}},\gamma\m {\tilde{P}}}'(\beta_\gamma^*)=0$ (recall~\eqref{eq:dl}), which shows that $\beta_\gamma^*$ is the maximizer of $\ell_{\gamma^{-1/2}\vect{\tilde{c}},\gamma\m {\tilde{P}}}(\beta)$. Therefore, Theorem~\ref{th:ellipsidInclusionUnitBall} ensures that $\ellipsoid{\vect c}{\m P}\subseteqeq \ellipsoid{\vect c_0}{{\gamma^{-1}} \m P_0}$, which is the inclusion~\eqref{eq:consequence_inclusion_1} in the statement. As a consequence, $\eqref{eq:consequence_inclusion_2}$ also holds as the ellipsoids therein are the same as those in \eqref{eq:consequence_inclusion_1} after a translation of $-\vect{c}_0$, a uniform scaling of $(-\ell_{\tilde{\vect{c}},\tilde{\m{P}}}^*)^{-1/2}$, and a translation of $\vect{c}_0$.\hfill\hfill\qed
\end{pf}

Notice that, different interpretations can be given to the inclusion~\eqref{eq:consequence_inclusion_1} in Corollary~\ref{cor:touching}: 
\begin{itemize}[leftmargin=*]
    \item $\ell_{\tilde{\vect{c}},\tilde{\m{P}}}^*>-1$: The ellipsoid $\ellipsoid{\vect c_0}{\m P_0}$ can be compressed by, at most, a factor of $(-\ell_{\tilde{\vect{c}},\tilde{\m{P}}}^*)^{1/2}<1$ and will still contain $\ellipsoid{\vect c}{\m P}$;
    \item $\ell_{\tilde{\vect{c}},\tilde{\m{P}}}^*=-1$: Corollary~\ref{cor:touching} becomes trivial, as we have $\gamma=1$ and there is already a contact point between the boundaries of  $\ellipsoid{\vect c_0}{\m P_0}$ and $\ellipsoid{\vect c}{\m P}$;
    \item $\ell_{\tilde{\vect{c}},\tilde{\m{P}}}^*<-1$:  The ellipsoid $\ellipsoid{\vect c_0}{\m P_0}$ must be inflated by, at least, a factor of $(-\ell_{\tilde{\vect{c}},\tilde{\m{P}}}^*)^{1/2}>1$ to contain $\ellipsoid{\vect c}{\m P}$.
\end{itemize}
The same applies to \eqref{eq:consequence_inclusion_2} but considering inflation and compression of $\ellipsoid{\vect c}{\m P}$ with respect to the $\vect {c}_0$.

Intuitively, Corollary~\ref{cor:touching} provides a simple method for determining what is the ellipsoid of minimum volume centered in $\vect{c}_0$ and with a shape defined by the spectrum of $\m P_0$ that contains $\ellipsoid{\vect{c}}{\m P}$. Alternatively, we can also compute the least sub-level set of the quadratic function $\vect{x}\mapsto (\vect{x}-\vect{c}_0)^\top\m P_0(\vect{x}-\vect{c}_0)$ containing the $n$-ellipsoid $\ellipsoid{\vect{c}}{\m P}$. This interpretation will be explored in the following sections, with applications in control design problems.

\section{Numerical Experiments}\label{sec:numerical}

\subsection{Implementation details}
To solve the optimization problem in Line~\ref{line:maximization} of Algorithm~\ref{algo:inclusion}, we propose a bisection algorithm (Algorithm~\ref{algo:bisection}) with stopping criteria guaranteeing inclusion without contact point ($\mathcal{E}\subset\Int{\mathcal{E}_0}$) and non-inclusion ($\mathcal{E}\not\subseteq\mathcal{E}_0$). Note that the inclusion with contact points ($\mathcal{E}\subseteqeq\mathcal{E}_0$) cannot be decided numerically, motivating the definition of $\mathcal{E}\subseteqeq^\epsilon\mathcal{E}_0$, which means that, given a machine precision $\epsilon$ of the computer, the algorithm cannot determine if $\mathcal{E}\subset\Int{\mathcal{E}_0}$ or $\mathcal{E}\not\subseteq\mathcal{E}_0$.

\SetNlSty{textbf}{}{ }\IncMargin{2em} 
\begin{algorithm2e}[tb]
  \SetAlgoLined
   \nl $l_0,u_0\gets \lambda_{\min}(\m P)^{-1},1-\vect c^\top \vect c$\;
   \nl \uIf{$\ell_{\vect c,\m P}(l_0) > -1$ {\normalfont \textbf{or}} $\ell_{\vect c,\m P}(u_0) > -1 $}{
   \nl \Return{$\mathcal{E}\subset \Int{\mathcal{E}_0}$}\;}
   \nl \uIf{$\ell_{\vect c,\m P}'(u_0) > 0$ }{
   \nl \Return{$\mathcal{E}\not\subseteq \mathcal{E}_0$}\;}\label{line:bisec-out-of-bounds}
   \nl $k\gets 0$\;
   \nl \While{$u_k-l_k>\epsilon$}{
   \nl $\beta_k = \tfrac{l_k+u_k}{2}$\;
   \nl \uIf{$\ell_{\vect c,\m P}(\beta_k) >-1$ \label{line:ifincluded}}{
   \nl \Return{$\mathcal{E}\subset\Int{\mathcal{E}_0}$}\;}
   \nl \uElseIf{$\ell_{\vect c,\m P}(\beta_k) < -1 + \ell_{\vect c,\m P}''(l_k)\tfrac{(u_k-l_k)^2}{2}$ \label{line:ifnotincluded}}{\nl \Return{$\mathcal{E}\not\subseteq \mathcal{E}_0$}\;}
   \nl \uIf{$\ell_{\vect c,\m P}'(\beta_k)<0$} 
   {\nl $l_{k+1},u_{k+1}\gets \beta_k, u_k$\;\label{line:ifnegatif}}
    \nl \uIf{$\ell_{\vect c,\m P}'(\beta_k)>0$}
   {\nl $l_{k+1},u_{k+1}\gets l_k, \beta_k$\;\label{line:ifpositif}}
   \nl $k\gets k+1$\;
   }
    \nl \Return{$\mathcal{E}\subseteqeq^\epsilon\mathcal{E}_0$}\;
  \caption{
    Test the inclusion of an ellipsoid $\mathcal{E} = \ellipsoid{\vect c}{\m P}$ in the ball $\mathcal{E}_0 = \ball{\vect 0}{1}$. 
    Returns either $\mathcal{E}\subset\Int{\mathcal{E}_0}$, $\mathcal{E}\not\subseteq\mathcal{E}_0$, or $\mathcal{E}\subseteqeq^\epsilon\mathcal{E}_0$.
  }
  \label{algo:bisection} 
\end{algorithm2e}\DecMargin{2em}
\begin{proposition}
Algorithm~\ref{algo:bisection} is correct and terminates in a finite number of steps.
\end{proposition}
\begin{pf}
Let $\beta^* = \argmax_{\beta\in \mathcal{I}_{\m P}}\ell_{\vect c,\m P}(\beta)$. 
Whenever $\beta^*\notin \mathcal{I}_{\vect c, \m P}$ we have $ \ell_{\vect c,\m P}'(1-\vect c^\top\vect c)>0$ and $\ell_{\vect c,\m P}(1-\vect c^\top\vect c)<-1$ implies, by Proposition~\ref{prop:upper_bound}, that $\mathcal{E}\not \subseteq\mathcal{E}_0$, which the algorithm deals with in line~\ref{line:bisec-out-of-bounds}. 
Let us consider the case $\beta^*\in \mathcal{I}_{\vect c,\m P}$.
In an arbitrary interval $[l,u]\subset\mathcal{I}_{\vect c,\m P}$ containing $\beta^*$ for which $l>1/\lambda_{\min}(\m P)$, the function $\ell_{\vect c,\m P}$ is locally $L$-smooth with $L = \max_{\beta\in [l,u]}|\ell_{\vect c,\m P}''(\beta)| = -\ell_{\vect c,\m P}''(l)$
since $\ell_{\vect c,\m P}''$ is negative and increasing in $[l,u]$.
Therefore, due to \cite[Lemma 1.2.3]{nesterov2018lectures} we have
$$\ell_{\vect c,\m P}(\beta^*)- \tfrac{L}{2}(\beta-\beta^*)^2 \le \ell_{\vect c,\m P}(\beta)\le \ell_{\vect c,\m P}(\beta^*)\quad \forall \beta \in [l,u]$$
and since $\beta^*\in [l,u]$, we obtain the following lower and upper bounds on $\ell^*_{\vect c,\m P}$ for all $\beta \in [l,u]$:
\begin{equation}\label{eq:bounds}
    \ell_{\vect c,\m P}(\beta)  \le \ell_{\vect c,\m P}(\beta^*)\le \ell_{\vect c,\m P}(\beta)- \tfrac{\ell_{\vect c,\m P}''(l)}{2}(u-l)^2.
\end{equation}

At each iteration $k$, lines~\ref{line:ifnegatif} or \ref{line:ifpositif} ensure that $\beta^*\in [l_k,u_k]$. 

If at a given iteration $k$, either $\ell_{\vect c,\m P}(\beta_k)\geq-1$, or $\ell_{\vect c,\m P}(\beta_k)- {\ell_{\vect c,\m P}''(l_k)}(u_k-l_k)^2/2<-1$, one can respectively conclude by \eqref{eq:bounds} that $\ell_{\vect c,\m P}^*>-1$ or $\ell_{\vect c,\m P}^*<-1$.

Whenever $\ell_{\vect c,\m P}^*>-1$, the continuity of $\ell_{\vect c,\m P}(\beta)$ ensures that the condition in line~\ref{line:ifincluded} will be satisfied at a given iteration, given that the interval $[l_k,u_k]$ converges to $\beta^*$ as $k\rightarrow\infty$. In the opposite case, when $\ell_{\vect c,\m P}^*<-1$, Algorithm~\ref{algo:bisection} also stops by satisfying the condition in line~\ref{line:ifnotincluded}. This holds by the same argument that $[l_k,u_k]$ converges to $\beta^*$, which also implies that $(u_k-l_k)^2\rightarrow 0$ and $\ell_{\vect c,\m P}''(l_k)\rightarrow\ell_{\vect c,\m P}''(\beta^*)<0$ as $k\rightarrow \infty$, concluding the proof. Finally,  if $\ell_{\vect c,\m P}^*=1$, the stop criterion of the main loop (i.e., $u_k-l_k\le\epsilon$, for a given precision $\epsilon>0$) will be satisfied and the algorithm returns $\mathcal{E}\subseteqeq^\epsilon\mathcal{E}_0$, implying that for this precision $\epsilon$ the inclusion with contact points may hold.\hfill\hfill\qed

\end{pf}

Finally, note that evaluating $\ell_{\vect c,\m P}(\beta)$ and its derivatives are computationally inexpensive operations. Once $\bar{c}_i^2,~i\in\support{\vect{\bar{c}}},$ and the eigenvalues of $\m P$ have been computed, the exact number of FLOPs to evaluate the function $\ell_{\vect c,\m P}$ for a given $\beta$ is $5|\support{\vect{\bar{c}}}|$, which, in the worst case, represents $5n$ FLOPs. 
Similarly, by performing some preliminary computations, the evaluation of $\ell'_{\vect c,\m P}$ and $\ell''_{\vect c,\m P}$ at a given $\beta$ requires respectively $5|\support{\vect{\bar{c}}}|$ and $6|\support{\vect{\bar{c}}}|$ FLOPs. Note that the evaluation of $\ell_{\vect c,\m P}$ and its derivative share several algebraic operations, which can be used to reduce the total number of computations.

\subsection{Performances}

Let us now compare the performance of  Algorithm~\ref{algo:inclusion} with those of solving the LMI condition~\eqref{eq:A_matrix} provided by \cite{boyd1994linear} with two conventional SDP solvers, namely, SDPA~\citep{yamashita2010high} and Mosek~\citep{mosek}. The underlying SDP problem tries to find $\beta\ge0$ such that the matrix $\m F(\beta)\succeq \m 0$ and, although a single variable is being searched for, the problem deals with an SDP restriction of size $n+1$.

\begin{figure}
    \centering
    \includegraphics[width=0.9\linewidth]{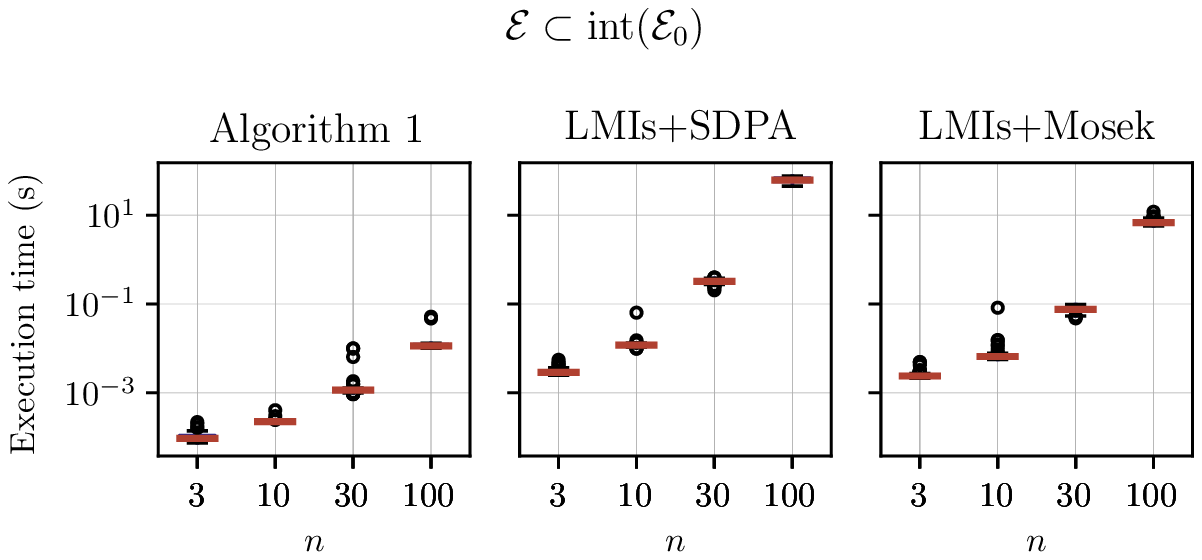}
    \includegraphics[width=0.9\linewidth]{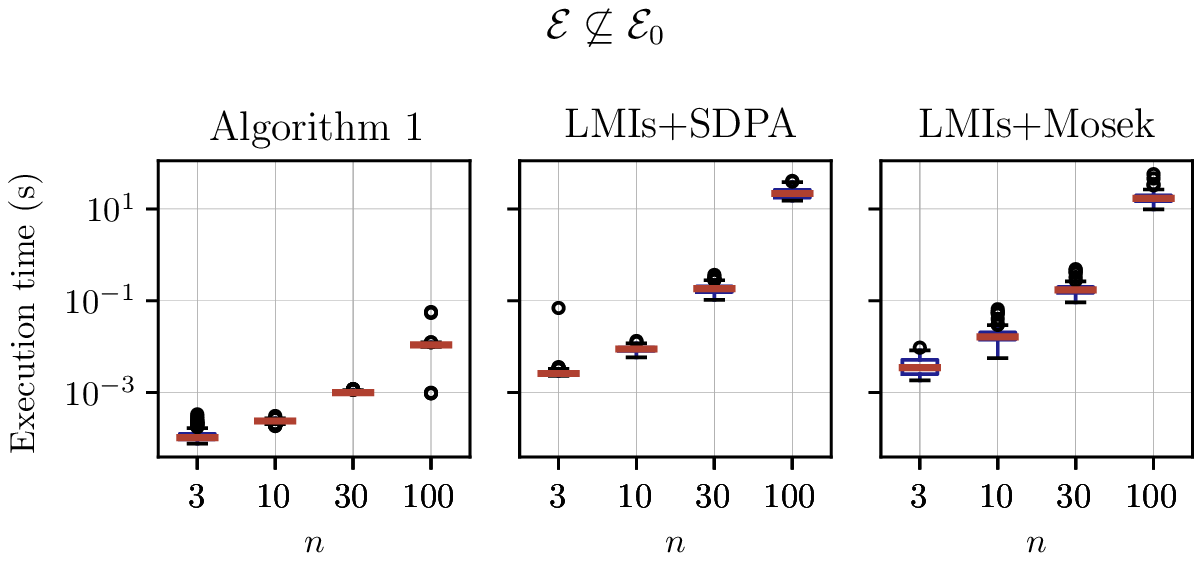}
    \caption{
    Comparison of the execution time for testing ellipsoid inclusion $\ellipsoid{\vect c}{\m P} \subseteq \ellipsoid{\vect c_0}{\m P_0}$ as a function of ellipsoid dimension $n$.}
    \label{fig:comparisonInclusion}
\end{figure}

Fig.~\ref{fig:comparisonInclusion} shows the average execution times (in seconds) for 200 randomly generated problems, each  consisting of two ellipsoids $\mathcal{E}=\ellipsoid{\vect{c}}{\m P}$ and  $\mathcal{E}_0=\ellipsoid{\vect{c}_0}{\m P_0}$. These were generated in a way that, in 100 cases, we have $\mathcal{E}\subset\Int{\mathcal{E}_0}$, and $\mathcal{E}\not\subseteq\mathcal{E}_0$ in the remaining ones.
Also, we made sure that the conditions in Proposition~\ref{prop:pretests} do not hold, so running the bisection Algorithm~\ref{algo:bisection} was required each time.
The comparison with the SDP solvers was repeated for $n$-ellipsoids with $n\in\{3,10,30,100\}$ to evaluate how these approaches scale up. These were performed in an Intel(R) Xeon(R) W-2295 CPU @ 3.00GHz with Julia  1.7.3 and using the optimization toolbox JuMP~\citep{DunningHuchetteLubin2017}. We computed that, in this benchmark, the execution times for Algorithm~\ref{algo:inclusion} were lesser than those for SDPA and Mosek in all cases. On average, our algorithm performs 27,  49,  162, and  2294 times faster than the competitors for ellipsoids of dimensions $3$, $10$, $30$, and $100$, respectively. In terms of memory consumption, Table~\ref{tab:mem} shows the average allocated memory (in kilobytes) of each case. This also demonstrates that our method is not only faster but also requires considerably fewer resources than the LMI-based approach.

\begin{table}
    \centering
    \caption{Memory Allocated on Average (in kB)}
    \!\!\!\begin{tabular}{c|c|c|c|c|c|c}
            &\multicolumn{3}{c|}{$\mathcal{E}\subset {\rm int}\mathcal{E}_0$} &\multicolumn{3}{c}{$\mathcal{E}\not\subseteq \mathcal{E}_0$} \\
        $n$ &  \!Alg.~\ref{algo:inclusion} \! &  SDPA   &  MOSEK  & \!\! Alg.~\ref{algo:inclusion} \! &  SDPA   &  MOSEK  \\
        3   &   8.2 &   192.6 &   153.5 &  13.8 &   203.3 &   152.1 \\
        10  &  26.5 &   733.2 &   577.8 &  36.4 &   733.1 &   576.0 \\
        30  & 118.4 &  4866.6 &  3792.1 & 138.6 &  4866.6 &  3790.3 \\
        100 & 865.8 & 51822.6 & 39608.1 & 901.8 & 51822.6 & 39606.2 \\
    \end{tabular}
    \label{tab:mem}
\end{table}
\subsection{Applications in Control Theory}\label{sec:application}
This section illustrates one possible application in the control theory of linear time-invariant (LTI) systems, namely, the computation of a control forward-invariant set for additive disturbances.  This problem has been extensively studied in the literature (e.g., see~\citep[Section~4.1]{blanchini1999set}) and, therefore, the goal of this section is not to provide a new method to tackle it  but rather to demonstrate the results from Corollary~\ref{cor:touching}.

For that, consider an LTI system
\begin{equation}
	\dot{\vect x} = \m A \vect x + \m B \vect u + \m H \vect w
\end{equation}
where $\vect x(t)\in\R^2$ is the state variable, and the signals $\vect u(t),~\vect w(t)\in\R$ are the control and the additive disturbance. This system is defined by the matrices
\begin{equation}
	\m A = \begin{pmatrix}
		0  & 1\\
		0.1 & 0.3
	\end{pmatrix},~ \m B =\begin{pmatrix}
	0\\0.5
\end{pmatrix} ,~ \m H =\begin{pmatrix}
-0.3\\0.6
\end{pmatrix} 
\end{equation}
and is controlled by a state-feedback LQR controller $\vect{u}(t) = \m K \vect{x}(t)$ synthesized for LQR parameters $\m Q = \m I$ and $ R = 1$, see \cite[p.~ 114]{boyd1994linear} for details on the LQR problem. The corresponding solution of the Algebraic Riccati Equation and feedback gain are 
\begin{equation}
	\m P = \begin{pmatrix}
		36.10& 42.36\\
		42.36& 72.98
	\end{pmatrix},~\m K = \begin{pmatrix}
	4.24& 7.30
\end{pmatrix},
\end{equation}
which allows us to define the closed loop matrix $\m A_{c} = \m A - \m B \m K$ and a Lyapunov function $v(\vect x) = \vect x^\top\m P \vect x$. 

\textbf{Bounded additive disturbances:} Due to the additive disturbance, the descent condition for this Lyapunov function does not hold everywhere. For an arbitrary $\vect w\in\mathcal{W}\subset \R$, this can be verified as
\begin{align}
	\dot v(\vect x,\vect w) &\coloneqq 2\vect{x}^\top (\m P \m A_{c})\vect{x} + 2\vect{x}^\top \m P \m H \vect{w}\nonumber\\
	&= -(\vect{x}-\m G\vect w)^\top \m S(\vect{x}-\m G\vect w) + r(\vect w)\label{eq:lyap_violated}
\end{align}
where $\m S = -\m A_{c}^\top \m P - \m  P\m A_c $,  $\m G = \m S^{-1}\m P\m H$ and $r(\vect w) = \vect w^{\top}\m G^{\top} \m S \m G\vect w$. From~\eqref{eq:lyap_violated}, one can conclude that $\dot{v}(\vect x,\vect w)\geq 0$ if and only if $\vect x \in \ellipsoid{\m G \vect w}{r(\vect w)^{-1}\m S}$. Due to the linearity of $\dot{v}(\vect x,\vect w)$ with respect to $\vect w$, when the disturbance $\vect w(t)$ takes values from a polytope $\mathcal{W}={\rm co} \{\vect w_1,\dots, \vect w_N\}$, the problem of finding the smallest sublevel set of $v(\vect{x})$ that is control forward-invariant reduces to finding the smallest $\gamma\geq 0$ such that 
\[\mathcal{B}_i\coloneqq\ellipsoid{\m G \vect w_i}{r(\vect w_i)\m S} \subseteq \mathcal{V}\coloneqq\ellipsoid{\vect 0}{\gamma^{-1} \m P}
\] for all $i=1,\dots, N$. In this context, Corollary~\ref{cor:touching} allows us to calculate $\gamma_i$ such that $\mathcal{B}_i \subseteqeq \mathcal{V}_i\coloneqq\ellipsoid{\vect 0}{\gamma_i^{-1} \m P}$ by comparison with the inclusion given in equation~\eqref{eq:consequence_inclusion_1}. Hence, we have that $\gamma = \max(\gamma_1,\dots,\gamma_N)$ is ensured to be the smallest such that $\mathcal{B}_i\subseteq\mathcal{V}$ for all $i=1,\dots, N$. Considering $\mathcal{W}= [-0.5,~0.5]$, we obtain $\gamma=1.137$. Notice that, for obtaining all $\gamma_1,\dots, \gamma_N$, the Cholesky decomposition and the spectral decomposition associated with the definition of the function~\eqref{eq:function_l} and the variable transformation~\eqref{eq:tildes} can be computed only once, given that the matrices defining the shape of all $\mathcal{B}_i$ are scalar multiples of $\m S$. Moreover, the symmetry of $\mathcal{V}$ ensures that  $\ellipsoid{\m G \vect w_i}{r(\vect w_i)^{-1}\m S} \subseteq \mathcal{V}$ implies that $\ellipsoid{-\m G \vect w_i}{r(\vect w_i)^{-1}\m S} \subseteq \mathcal{V}$, which helps to reduce the number of executions of the bisection algorithm to calculate $\ell_{\vect c, \m P}^*$ in Corollary~\ref{cor:touching}.

\begin{figure}
\centering
	\includegraphics[width=0.77\linewidth]{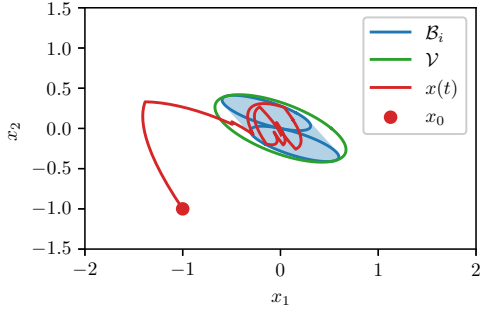}
	\caption{For an disturbed LTI system presented in Section~\ref{sec:application}, we computed with the results from Corollary~\ref{cor:touching} the smallest forward-invariant level set $\mathcal{V}$ of the Lyapunov function $v(\vect x)$.} \label{fig:control}
\end{figure}

 In Figure~\ref{fig:control}, the ellipsoidal sets defined in this section are illustrated, along with a trajectory undergoing a random disturbance and starting at $\vect{x}_0 = [-1~~-1]^\top$. The light blue shaded area represents the set ${\rm co}\{\mathcal{B}_1, \mathcal{B}_2\}$, where the decreasing property of the Lyapunov function fails for some $\vect{w}\in\mathcal{W}$.
 
In summary, Corollary~\ref{cor:touching} allows the verification that the Lyapunov function $v(\vect x)=\vect{x}^\top \m P \vect{x}$ strictly decreases in  $\R^2\setminus\mathcal{V}$ despite the persistent disturbance. This verification is done efficiently by performing the following numerical operations: one Cholesky decomposition, one spectral decomposition, and one bisection algorithm for maximizing a concave scalar function in a compact interval.
\section{Conclusion and Future Work}
We presented a new method to verify the inclusion of $n$-ellipsoids, which consists in the maximization of a scalar concave and smooth function~\eqref{eq:function_l}. This function and its derivatives can be computed in $\mathcal{O}(n)$ floating-point operations and the interval~\eqref{eq:interval_I} where its maximum lies is a subset of $[0,1]$. Therefore, we proposed a bisection-based algorithm (Algorithm~\ref{algo:bisection}) allowing us to decide whether the inclusion holds. A benchmark with methods based on LMI constraints tackled by two off-the-shelf SDP solvers is carried out, showing that we outperform the LMI-based approach. We also present an application in the field of control theory.

The source codes for the numerical experiments carried out in this paper are available in the following repository: \texttt{\small https://github.com/egidioln/EllipsoidInclusion.jl}.

For future work, we plan to generalize this approach to verify the emptiness of the intersection of ellipsoids and also other quadrics. We also seek to apply these results to model and data-driven abstraction-based control.

\bibliography{main}

\appendix

\end{document}